\documentclass [10pt]{article}
\usepackage{amsmath}

\usepackage{amsfonts, amssymb}

\usepackage[cp1251]{inputenc} 

\def\be#1\ee{\begin{equation}#1\end{equation}}

\newtheorem{thm}{Theorem}
\newtheorem{lem}[thm]{Lemma}

\def\P{{\mathbb{P}}}
\def\R{\mathbb{R}}
\def\E{\mathbb{E}\,}

\def\Z{{\mathbb Z}}

\def\AAA{{\mathcal{E}}}

\newcommand{\eps}{\varepsilon}

\def\HH{{\mathcal H}}

\def\tX{{\widetilde{X}}}
\def\tY{{\widetilde{Y}}}

\def \=L{\ {\buildrel\hbox{\scriptsize d }\over =}\ }

\begin{document}

\title{ On the distribution of the last exit time over a slowly growing linear boundary for a Gaussian process\footnote{The work
of M.A. Lifshits was supported by RFBR-DFG grant 20-51-12004.}
}

\author{
N.A.\,Karagodin\footnote{St.Petersburg State University, Department of Mathematics and Computer Sciences,
199034, St.Petersburg, Universitetskaya emb., 7–9. {\tt email: nikitus20@gmail.com}.}
,
M.A.Lifshits\footnote{St.Petersburg State University, Department of Mathematics and Computer Sciences,
199034, St.Petersburg, Universitetskaya emb., 7–9.  {\tt email: mikhail@lifshits.org}.}
}

\date{}

\maketitle

\begin{abstract} For a class of Gaussian stationary processes, we prove a limit theorem on the convergence
of the distributions of the scaled last exit time over a slowly growing linear boundary.
The limit is a double exponential (Gumbel) distribution.
\end{abstract}


\medskip
\bigskip

\bigskip

\section{Introduction}

\subsection{Problem setting}

Consider a stationary Gaussian process with continuous trajectories and its  "last exit time over a linear boundary",
i.e. the last instant when the process hits a line $at$, where $t$ denotes time  and $a>0$ is a drift (or trend) parameter.
After this instant, the process stays forever {\it under} the line. We are interested in the asymptotic distribution of the last exit time
when the trend  $a$ goes to zero. In this work, we will prove a limit theorem on the convergence of the distribution of the properly
centered and scaled last exit time to a double exponential (Gumbel) law.

A special case of this problem, for a particular process, emerged in recent works \cite{ABL19,ABL20} providing a mathematical study
of a physical model (Brownian chain break). Quite naturally, a question is raised, whether it is possible to extend that result to a sufficiently wide class
of processes. This is what we do here.

As far as we know, the problem setting handling the small trend is new, although the last exit time is a sufficiently popular object in the problems of
economical applications such as studies of ruin probabilities. In those settings, however, as a rule, one considers processes with stationary
increments and trend is fixed, see \cite{Deb,Hus,ParRab}.

\subsection{Main result}

Let $Y(t),t\in \R$, be a real-valued centered stationary Gaussian process with covariance function $\rho(t):=\E [Y(t)Y(0)]$.
We make two assumptions on the covariance function: at zero
\begin{eqnarray} \label{assumption_0}
    \rho(t)=v^2(1-Q|t|^\alpha+o(|t|^\alpha)), \quad \textrm{as } t\to 0,
\end{eqnarray}
for some $v>0, Q>0$, $\alpha\in(0,2]$,
and at infinity
\begin{eqnarray} \label{assumption_infty}
  \rho(t)= o( (\ln t)^{-1} ),   \quad \textrm{as } t\to +\infty.
\end{eqnarray}

Recall that relation \eqref{assumption_0} appears in the following lemma that will serve as one of the two basic tools
in our calculations.

\begin{lem} \label{lem:pikpit} (Pickands--Piterbarg lemma).
Let $Y(t),t\in \R$, be a real-valued centered stationary Gaussian process satisfying conditions
\eqref{assumption_0} and
\[
   \limsup_{t\to\infty} \rho(t) <1.
\]
Then
\[
  \P \left\{ \max_{s\in[0,t]} Y(s) \ge x \right\}
  \sim \frac{Q^{1/\alpha}\HH_\alpha}{\sqrt{2\pi}} \cdot t \cdot (x/v)^{2/\alpha-1} e^{-x^2/2v^2}
\]
for all $x$ and $t$ such that the right hand side tends to zero and $tx^{2/\alpha}\to\infty$. Here $\HH_\alpha$ are
Pickands constants (in particular, $\HH_1=1$,$\HH_2=\pi^{-1/2}$).
\end{lem}

Here and elsewhere throughout the paper  $f\sim g$ stands for $\lim\tfrac fg=1$.
A first version of this lemma with fixed $t$ was obtained by Pickands \cite{Pick}, while
this version with variable $t$ (which is very important for our goals) is due to Piterbarg,
see \cite[lecture 9]{Pit}.
\bigskip

Define the $Y$'s last exit time over a linear boundary as
\[
   T= T(\eps):=\max\{t: Y(t)=\eps\,t\}.
\]

The main result of our work is as follows.

\begin{thm} \label{t:main} Let $Y(t),t\in \R$, be a real-valued centered stationary Gaussian process satisfying
assumptions \eqref{assumption_0} and \eqref{assumption_infty}. Let
$c:=\frac{Q^{1/\alpha}\HH_\alpha}{\sqrt{2\pi}}$,  $\eps_v:=\frac{\eps}{v}$.
Then for each $r\in \R$ it is true that
\[
  \lim_{\eps\to 0} \P\left\{ \frac{T(\eps)-A_\eps}{B_\eps}\le r \right\} = \exp(-c\exp(-r)),
\]
with the scaling constants
\begin{eqnarray*}
   A_\eps&:=& \eps_v^{-1} \left( \sqrt{-2\ln\eps_v} + \left(\frac 1{\alpha}-1\right)\,
   \frac{\ln(-2\ln\eps_v)}{\sqrt{-2\ln\eps_v}} \right),
\\
   B_\eps &:=& \left(\eps_v\sqrt{-2\ln\eps_v}\right)^{-1}.
\end{eqnarray*}
\end{thm}
\medskip

We stress that the research technique of the initial work \cite{ABL19} was based on the exponential strong mixing property satisfied
by the process studied there. In our work, the assumption is only imposed on the covariance function (which is
by far easier to check than strong mixing) and the condition of the covariance decay is logarithmic  instead of the exponential
one.
\medskip

The double exponential law usually emerges in the studies of maxima of the identically distributed random variables.
Amazingly,  in our problem it is related to the maxima of {\it non}-identically distributed variables.

\section{Proofs}

\subsection{Proof of Theorem \ref{t:main}}

By the linear variable change $Y(t)= v{\widetilde Y}(Q^{1/\alpha}t)$ one may reduce the problem
to the case $v=Q=1$, $\eps_v=\eps$, which is considered in the following.

Let us fix $r\in \R$ and let
\[
   \tau=\tau(\eps,r):= A_\eps +B_\eps r.
\]
The theorem's statement is equivalent to
\[
  \lim_{\eps\to 0} \P\left\{ T(\eps)\le \tau \right\} = \exp(-c\exp(-r)).
\]
In order to prove this, we cover the halfline $[\tau,\infty)$ by the following system of sets:

\begin{itemize}
\item
the halfline $[\sigma,\infty)$, where $\sigma:=A_\eps +B_\eps R$ and $R=R(\eps)$ slowly tends to infinity.
The choice of $R$ will be further specified at the end of the proof.
\item
long intervals $L_i=[(\ell+s)i,(\ell+s)i+\ell]$, $i\in\Z$, of some length  $\ell$,
\item
shorter intervals $S_i=[(\ell+s)i+\ell,(\ell+s)(i+1)]$,$i\in\Z$, of length $s$.

\end{itemize}
In fact the main part of process exits over the linear boundary will occur on the long intervals,
 while the shorter intervals placed between the long ones play the role of separators
 providing weak dependence between the process values on different long intervals.

The interval lengths $\ell=\ell(\eps)$, $s=s(\eps)$ must satisfy the relations

\begin{eqnarray} \label{ls0}
   \ln s \sim |\ln(-\eps)|,
\\  \label{ls1}
    s/\ell\to 0,
\\  \label{ls2}
    \eps \, \sqrt{-\ln\eps} \, \ell \to 0.
\end{eqnarray}

The parameter $R$ should grow to infinity so slowly that
\begin{eqnarray} \label{sigma0}
	\tau \sim \sigma.
\end{eqnarray}	

Let
\[
   X_i^\eps:= \max_{t\in L_i} Y(t); \quad  V_i^\eps:= \max_{t\in S_i} Y(t).
\]
By using stationarity, we infer from Pickands--Piterbarg lemma the asymptotics
$\P\{ X_i^\eps\ge x\} \sim c\, \ell\, x^{2/\alpha-1} \exp(-x^2/2)$
and $\P\{ V_i^\eps\ge x\} \sim c\, s\, x^{2/\alpha-1} \exp(-x^2/2)$,
as soon as the corresponding right hand sides tend to zero.
Here $c=\tfrac{\HH_\alpha}{\sqrt{2\pi}}$.

Define the index sets
\begin{eqnarray*}
   I_1 &:=& \{i:\ \substack{(\ell+s)i+\ell\ge\tau \\(\ell+s)i < \sigma }\},
\\
   I_2 &:=& \{i:\ \substack{(\ell+s)i \ge\tau \\(\ell+s)i+\ell < \sigma}\},
\\
   I_3 &:=& \{ i: \substack{(\ell+s)(i+1) \ge\tau \\(\ell+s)i+\ell < \sigma}\},
\end{eqnarray*}
chosen so that the following inclusions hold:
\be \label{inclusions}
     \bigcup_{i\in I_2}  L_i  \subset  [\tau,\sigma]
     \subset  \big( \bigcup_{i\in I_1} L_i \big) \cup \big( \bigcup_{i\in I_3} S_i \big).
\ee
In the first inclusion one considers the long intervals belonging to $[\tau,\sigma]$;
in the second one, the long intervals and the short intervals separating them cover $[\tau,\sigma]$.

Let us define the events related to the exits of our process over the linear boundary:
\begin{eqnarray*}
     \AAA_1 &:=& \bigcup_{i\in I_1} \{ X_i^\eps \ge(\ell+s)i\eps \},  
\\
     \AAA_2 &:=&  \bigcup_{i\in I_2}  \{ X_i^\eps\ge (\ell+s)(i+1)\eps \},  
\\
     \AAA_3 &:=& \bigcup_{i\in I_3}   \{ V_i^\eps\ge (\ell+s)i\eps \},  
\\
     \AAA_4 &:=&  \{ \exists\, t>\sigma : Y(t) \ge \eps t \}.
\end{eqnarray*}

By using inclusions \eqref{inclusions} and monotonicity of the linear function,
it is easy to see that the following bounds are true:
\begin{eqnarray*}
   \P\{ T(\eps)>\tau \} &=& \P\{ \exists\, t> \tau: Y(t) \ge \eps t \} \le \P\{\AAA_1\}+\P\{\AAA_3\}+\P\{\AAA_4\},
\\
   \P\{ T(\eps)>\tau \} &\ge& \P\{\AAA_2\}.
\end{eqnarray*}
Therefore, it is sufficient to prove that, as $\eps\to 0$,  we have
\begin{eqnarray*}
  \P\{\AAA_1\}, \P\{\AAA_2\}  &\rightarrow& 1-\exp(-c\exp(-r)),
\\
  \P\{\AAA_3\}, \P\{\AAA_4\} &\rightarrow& 0.
\end{eqnarray*}

Let us first prove that the probabilities of the events $\AAA_1$ and $\AAA_2$ are almost equal,
thus it will be enough to find the limit of $\P\{\AAA_1\}$.
Indeed, let the indices $m$ and $n$ be such that $I_1=[m,n]$. Then
\begin{eqnarray*}
  \P\{\AAA_2\} &\le& \P\{\AAA_1\} = \P\left\{ \bigcup_{i=m}^n \{ X_i^\eps \ge(\ell+s)i\eps \} \right\}
\\
  &\le& \P\left\{X_m^\eps \ge(\ell+s)m\eps \right\} + \P\left\{X_{m+1}^\eps \ge(\ell+s)(m+1)\eps \right\}
\\
  &&  + \P\left\{ \bigcup_{i=m+2}^n \{ X_i^\eps \ge(\ell+s)i\eps \} \right\}
\\
   &\le& 2\, \P\left\{X_m^\eps \ge(\ell+s)m\eps \right\}
  + \P\left\{ \bigcup_{j=m+1}^{n-1} \{ X_{j+1}^\eps \ge(\ell+s)(j+1)\eps \} \right\}
\\
  &\le& 2\, \P\left\{X_m^\eps \ge(\ell+s)m\eps \right\} +   \P\{\AAA_2\},
\end{eqnarray*}
where in the penultimate inequality we used the stationarity of the sequence $X_i^\eps$ following
from the stationarity of the process $Y$.
For the remaining term we use Pickands--Piterbarg bound and obtain
\[
     \P\left\{X_m^\eps \ge(\ell+s)m\eps \right\} \le \P\left\{X_m^\eps \ge (\tau-\ell) \eps \right\}
     \sim
     c \ell [(\tau-\ell)\eps]^{2/\alpha-1} \exp\{ -[(\tau-\ell)\eps]^2/2\}.
\]
Let us take into account that the definitions of $\tau$ and $\sigma$ yield the following three relations:
\begin{eqnarray} \label{tauassymp_1}
   \tau &\sim& \eps^{-1} \sqrt{-2\ln\eps},
\\ \label{tauassymp_0}
   \frac{1}{\eps}(\tau\eps)^{2/\alpha-2}\exp\{-(\tau\eps)^2/2\} &\sim& e^{-r},
\\ \label{sigmaassymp_0}
   \frac{1}{\eps}(\sigma\eps)^{2/\alpha-2}\exp\{-(\sigma\eps)^2/2\} &=& o(1).
\end{eqnarray}
We stress that equation  \eqref{tauassymp_0} is a key for the choice
of the scaling constants $A_\eps$ and $B_\eps$ in the theorem assertion.

We use the bound \eqref{ls2} and obtain
$\ell/\tau\to 0$, $\ell\tau\eps^2\to 0$. Therefore,
\begin{eqnarray*}
   \ell [(\tau-\ell)\eps]^{2/\alpha-1} \exp\{-[(\tau-\ell)\eps]^2/2\}
   &\sim&
   \ell (\tau\eps)^{2/\alpha-1} \exp\{ -(\tau\eps)^2/2\}
 \\
   &\sim&
   \ell\tau\eps^2 e^{-r}\to 0.
\end{eqnarray*}

We conclude that $\P\left\{X_m^\eps \ge(\ell+s)m\eps \right\}\to 0$, so that the difference
between $\P\{\AAA_1\}$ and $\P\{\AAA_2\}$ is indeed negligible.
\medskip

In the sequel, we will many times use the following technical lemma. Its proof is postponed to
Section \ref{ss:proof_prop}.

\begin{lem} \label{l:prop}
	For each $\alpha\neq 0$ and all $\theta(\eps), a(\eps), b(\eps)$ such that, as $\eps \to 0$, one has
	$\theta\eps \to \infty$, $a = o(\theta)$, $\theta a \eps^2\to 0$, it is true that
	\[
      	\sum_{i: ai+b \ge\theta }^{\infty} [(ai+b)\eps]^{2/\alpha-1} \exp \{-[(ai+b)\eps]^2/2\}
	     \sim \frac{1}{a \eps}(\theta \eps)^{2/\alpha-2}\exp\{-(\theta\eps )^2/2\}.
	\]
\end{lem}

Let us evaluate $\P\{\AAA_3\}$. By Pickands--Piterbarg asymptotics, we have
\begin{eqnarray*}
   \P\{\AAA_3\} &\le& \sum_{i\in I_3 } \P\{ V_i^\eps\ge (\ell+s)i\eps \}
\\
   &\le& c\, s\, \sum_{i: (\ell+s)(i+1)\geq \tau }    [(\ell+s)i\eps]^{2/\alpha-1}\exp \{- [(\ell+s)i\eps]^2/2\} (1+o(1)).
\end{eqnarray*}

In order to find the asymptotic behavior of this sum, we apply Lemma \ref{l:prop} with parameters
$a=\ell+s, b=0, \theta = \tau-\ell-s$. Then, by using \eqref{ls1}, \eqref{ls2}, and \eqref{tauassymp_1}
we have $a \sim \ell, \theta \sim \tau, (\theta \eps)^2 = (\tau\eps)^2+o(1)$.
Therefore, Lemma \ref{l:prop}, relation \eqref{tauassymp_0},  and assumption \eqref{ls1} yield
\[
    \frac{cs}{a\eps}(\theta\eps)^{2/\alpha-2}\exp\{-(\theta\eps)^2/2\}
    \sim
    \frac{cs}{\ell\eps}(\tau\eps)^{2/\alpha-2}\exp\{-(\tau\eps)^2/2\}
    = o(1).
\]
\medskip

Evaluation of $\P\{\AAA_4\}$ goes along the same lines. By splitting the halfline $[\sigma,\infty)$
into the intervals of the unit length, we obtain
\begin{eqnarray*}
  && \P\{\AAA_4\} \le \sum_{j=0}^\infty  \P\{  \max_{t\in [ \sigma+ j, \sigma+j+1]} Y(t)> \eps (\sigma + j) \}
\\
   &\le&  c\, \sum_{j=0}^\infty  [(\sigma+ j)\eps]^{2/\alpha-1}\exp \{- [\sigma + j)\eps]^2/2\} (1+o(1)).
\end{eqnarray*}

We use Lemma \ref{l:prop} with parameters $a=1, b=\sigma, \theta = \sigma$. By applying \eqref{sigmaassymp_0},
we obtain the asymptotics
\[
   \frac{c}{\eps}(\sigma\eps)^{2/\alpha-2}\exp\{-(\sigma\eps)^2/2\}
   = o(1).
\]

\bigskip

The subsequent estimates use  the effect of weak dependence of the values of the process $Y$ at distant times.
Our main tool here is the following classical inequality due to Slepian (see e.g.,
\cite[\S 14]{Lif_GSF}, \cite[lecture 2]{Pit}).

\begin{lem} Let $(U_1,...,U_n)$ and  $(V_1,...,V_n)$ be two centered Gaussian vectors such that
$\E U_j^2=\E V_j^2$, $1\le j\le n$, and $\E (U_i U_j)\le \E (V_i V_j)$, $1\le i,j\le n$.
Then for each $r\in \R$  one has
\[
    \P\left\{ \max_{1\le j\le n} U_j \ge r \right\}
    \ge  \P\left\{ \max_{1\le j\le n} V_j \ge r \right\}.
\]
\end{lem}

One may write this inequality in a slightly more general form (see \cite[lecture 2]{Pit}): under assumptions of
Slepian lemma, for all non-negative $r_1,...,r_n$ one has
\[
   \P\left\{ \exists j:\ U_j \ge r_j \right\}
    \ge  \P\left\{ \exists j:\ V_j \ge r_j \right\}.
\]
This fact follows by application of Slepian inequality to the vectors
$(\tfrac{U_1}{r_1},...,\tfrac{U_n}{r_n})$, $(\tfrac{V_1}{r_1},...,\tfrac{V_n}{r_n})$
and $r=1$.

The latter inequality obviously extends to the Gaussian processes with continuous trajectories defined
on a metric space (by the way, the processes satisfying assumption \eqref{assumption_0} belong
to this class). Namely, let $\{U(t),t\in T\}$ and $\{V(t),t\in T\}$
be two Gaussian processes with continuous trajectories defined on a common metric space $T$.
Let $\E U(t)^2=\E V(t)^2$, $t\in T$, and $\E (U(t_1) U(t_2))\le \E (V(t_1) V(t_2))$, $t_1,t_2\in T$.
Then for all compact sets $T_1,...,T_n$ in $T$ and for all non-negative $r_1,...,r_n$ it is true that
\be \label{slepian}
    \P\left\{ \bigcup_{j=1}^n  \left\{ \max_{t\in T_j} U(t) \ge r_j \right\} \right\}
    \ge   \P\left\{  \bigcup_{j=1}^n  \left\{\max_{t\in T_j} V(t) \ge r_j\right\} \right\}.
\ee
\bigskip

Now we may proceed to the proof of the remaining claim
\be
    1- \P\{\AAA_1\}  \rightarrow \exp(-c\exp(-r)), \qquad \eps\to 0.
\ee
We provide the corresponding upper and lower bounds. In both cases we will use
Slepian inequality in the form \eqref{slepian}.
\medskip

{\it Upper bound.}

Let us compare our process $Y$ with an auxiliary process $Z$ which is defined as follows.
First, let us consider a process $\widetilde{Y}(t), t \in \cup L_i$
which consists of independent copies of $Y(t)$ on the intervals $L_i$.

Further, let
\[
  \delta^2=\delta^2(\eps):= \sup_{t\ge s(\eps)} \left|\E[ Y(t) Y(0) ]\right|.
\]
Taking into account the correlation decay assumption \eqref{assumption_infty} and assumption \eqref{ls0} concerning
the choice of $s$, we have
\be \label{deltas}
   \delta^2= o((\ln s)^{-1}) =  o((- \ln \eps)^{-1}).
\ee

Let $\xi$ be an auxiliary standard normal random variable independent with the process $\widetilde{Y}$.
We define the centered Gaussian process $Z(t)$, $t \in \cup L_i$, by the equality
\[
   Z(t) := \sqrt{1-\delta^2} \widetilde{Y}(t) + \delta \xi.
\]
Then for all $t$ the variances are equal: $\E Y(t)^2= \E Z(t)^2=1$.
For covariances we have the following inequalities:
\begin{itemize}
\item
    for $t_1$ and $t_2$ that belong to the same interval $L_i$ we have
	\begin{eqnarray*}
		\E[Z(t_1)Z(t_2)]&=&
		\E\left[\left(\sqrt{1-\delta^2} \tY(t_1) + \delta \xi\right) \left(\sqrt{1-\delta^2} \tY(t_2)+\delta \xi\right)\right]
		\\&=&(1-\delta^2)\E[Y(t_1)Y(t_2)]+\delta^2\ge \E[Y(t_1)Y(t_2)],
	\end{eqnarray*}
	where the last inequality follows from $\E [Y(t_1)Y(t_2)] \le \sqrt{\E Y(t_1)^2 \E Y(t_2)^2}$ $=1$,
\item
    for $t_1$ and $t_2$ that belong to different intervals $L_i$ and $L_j$,  by the definition of
    $\delta$ and by intervals' construction we have
	\begin{eqnarray*}
		\E[Z(t_1)Z(t_2)]&=&
		\E\left[\left(\sqrt{1-\delta^2} \tY(t_1) + \delta \xi\right) \left(\sqrt{1-\delta^2} \tY(t_2)+\delta \xi\right)\right]
		\\&=& \delta^2\ge \E[Y(t_1)Y(t_2)].
	\end{eqnarray*}
\end{itemize}

Let $\tX_i^\eps := \max_{t\in L_i} \tY(t)$. By applying Slepian inequality \eqref{slepian}
to the processes $Y$ and $Z$, we obtain
\[
   \P\{\AAA_1\}
   = \P\left\{ \bigcup_{i\in I_1} \{ X_i^\eps \ge(\ell+s)i\eps \} \right\}
   \ge \P\left\{ \bigcup_{i\in I_1} \{ \sqrt{1-\delta^2}\tX_i^\eps + \delta \xi \ge(\ell+s)i\eps \} \right\}.
\]
Let us pass to the complementary events; for every $h=h(\eps)>0$ the following elementary bound holds,
\begin{eqnarray} \nonumber
	1-\P\{\AAA_1\}
    &=& \P\left\{\bigcap_{i\in I_1} \{X_i^\eps\le (\ell+s)i\eps\}\right\}
    \\ \nonumber
    &\le& \P\left\{\bigcap_{i\in I_1} \{\sqrt{1-\delta^2} \tX_i^{\eps}+\delta \xi \le (\ell+s)i\eps\}\right\}
	\\ \nonumber
    &\le& \P\left\{\bigcap_{i\in I_1} \{\sqrt{1-\delta^2} \tX_i^{\eps}\le (\ell+s) i \eps + h\eps\}\right\}+\P\{\delta \xi \le -h\eps\}
	\\ \label{hepsdelta_prob}
    &=& \prod_{i\in I_1} \P\left\{ X_i^\eps\le \frac{(\ell+s) i \eps+h\eps}{\sqrt{1-\delta^2}} \right\}+\P\{\xi \le -h\eps/\delta\},
\end{eqnarray}
where the last equality holds because $\tX_i^\eps$  are independent copies of $X_i^\eps$.
We choose the level $h=h(\delta,\eps)$ so that
\begin{eqnarray} \label{h_large}
   && h\eps/\delta \to \infty,
\\  \label{h_small}
   && h\eps \cdot \sqrt{-2\ln\eps} \to 0,
\end{eqnarray}
which is possible under \eqref{deltas}.
It also follows from \eqref{h_small}  that
\[
  h/\tau = (h\eps) / (\tau\eps) \sim (h\eps)/ \sqrt{-2\ln \eps} \to 0.
\]

Due to \eqref{h_large}, the last term in our bound  \eqref{hepsdelta_prob} is negligible.
We check now that the product converges to $\exp(-c\exp(-r))$.
Taking the logarithm and passing to the complementary events, we see that it is necessary to
prove the convergence
\[
     \sum_{i\in I_1} \P\left\{ X_i^\eps\ge \frac{(\ell+s) i \eps+h\eps}{\sqrt{1-\delta^2}} \right\} \to c\exp(-r).
\]
By Pickands--Piterbarg lemma this is equivalent to
\[
     c \,\ell\, \sum_{i\in I_1} \left( \frac{(\ell+s) i \eps+h\eps}{\sqrt{1-\delta^2}} \right)^{2/\alpha-1}
     \exp\left(-\left(\frac{(\ell+s) i \eps+h\eps}{\sqrt{1-\delta^2}} \right)^2/2\right).
\]
We represent this expression as a difference of two sums
\be \label{sum_1}
	c \,\ell\, \sum_{i: (\ell+s)i+\ell\ge \tau} \left( \frac{(\ell+s) i \eps+h\eps}{\sqrt{1-\delta^2}} \right)^{2/\alpha-1}
	\exp\left(-\left(\frac{(\ell+s) i \eps+h\eps}{\sqrt{1-\delta^2}} \right)^2/2\right)
\ee
and
\be \label{sum_2}
c \,\ell\, \sum_{i: (\ell+s) i \ge \sigma } \left( \frac{(\ell+s) i \eps+h\eps}{\sqrt{1-\delta^2}} \right)^{2/\alpha-1}
\exp\left(-\left(\frac{(\ell+s) i \eps+h\eps}{\sqrt{1-\delta^2}} \right)^2/2\right)
\ee
The asymptotics of the first sum follows from Lemma \ref{l:prop} applied with parameters
$a=\frac{\ell+s}{\sqrt{1-\delta^2}}\sim \ell$ and $\theta=\frac{\tau - \ell +h}{\sqrt{1-\delta^2}}\sim \tau$,
where equations \eqref{h_small}, \eqref{tauassymp_1}, \eqref{ls2} and \eqref{deltas} yield
\begin{eqnarray*}
  (\theta\eps)^2&=& \frac{(\tau-\ell+h)^2\eps^2}{1-\delta^2} = \frac{\tau^2\eps^2}{1-\delta^2} + O(\tau (h+\ell) \eps^2)
\\
  &=&  \tau^2\eps^2 +  \frac{\tau^2\eps^2\delta^2}{1-\delta^2} + O(\tau (h+\ell) \eps^2)
  = \tau^2\eps^2 + o(1).
\end{eqnarray*}
By this relation and \eqref{tauassymp_0}, Lemma \ref{l:prop}
provides the following asymptotics for \eqref{sum_1}:
\[
    c \,\ell\, \frac{1}{a \eps}(\theta \eps)^{2/\alpha-2}\exp\{-(\theta\eps )^2/2\}
    \sim  c \,\ell\, \frac{1}{\ell \eps} (\tau \eps)^{2/\alpha-2}\exp\{-(\tau\eps )^2/2\}
    \sim ce^{-r}.
\]
Similarly, Lemma \ref{l:prop} applied with parameters $a=\frac{\ell+s}{\sqrt{1-\delta^2}}\sim \ell$
and $\theta = \frac{\sigma+h}{\sqrt{1-\delta^2}}\sim \sigma$ provides an asymptotics for \eqref{sum_2}.
Since $(\theta\eps)^2 = (\sigma\eps)^2+o(1)$,  by using \eqref{sigmaassymp_0} we obtain
\[
    c\ell\frac{1}{a\eps}(\theta\eps)^{2/\alpha-2}\exp\{-(\theta\eps)^2/2\}
    \sim c\ell\frac{1}{\ell\eps}(\sigma\eps)^{2/\alpha-2}\exp\{-(\sigma\eps)^2/2\}=o(1).
\]
Substraction of sums' asymptotics implies the required upper bound for $1-\P\{\AAA_1\}$.
\medskip

{\it Lower bound.}

In order to obtain an opposite bound for $1-\P\{\AAA_1\}$, we will introduce and compare two more auxiliary processes
$Y_1$, $\tY_1$.
Let $\xi$ be an auxiliary standard normal random variable independent with the process $Y$.
Let $Y_1(t) := Y(t)+\delta\xi$, $t\in\cup L_i$.  Furthermore, let us consider a sequence of independent standard Gaussian
random variables
$\xi_i$ independent of $\tY(t)$  and let
\[
    \tY_1(t) := \tY(t)+\delta\xi_i, \qquad t\in L_i.
\]
Then for all $t$ we have the equality of variances: $\E Y_1(t)^2= \E\tY_1(t)= 1+\delta^2$.
For covariances we have the following inequalities:
\begin{itemize}
	\item for $t_1$ and $t_2$ that belong to the same interval $L_i$ we have
\[
	\E[Y_1(t_1)Y_1(t_2)] = \E[\tY_1(t_1)\tY_1(t_2)],
\]
	\item for $t_1$ and $t_2$ that belong to different intervals $L_i$ and $L_j$ we have
\begin{eqnarray*}
      \E[Y_1(t_1)Y_1(t_2)] &=& \E[(Y(t_1)+\delta\xi)(Y(t_2)+\delta\xi)]=\E[Y(t_1)Y(t_2)]+\delta^2
\\
      &\ge& 0=	\E[\tY_1(t_1)\tY_1(t_2)].
\end{eqnarray*}
\end{itemize}

We choose $h=h(\delta,\eps)$ as before, i.e. satisfying assumptions
\eqref{h_large} and \eqref{h_small}.

Slepian inequality \eqref{slepian} yields
\begin{eqnarray*}
 &&  \P \left\{ \bigcup_{i\in I_1} \{ \tX_i^\eps+\delta \xi_i \ge (\ell+s) i \eps-h \eps \} \right\}
 = \P \left\{ \bigcup_{i\in I_1} \{ \max_{t\in L_i} \tY_1(t) \ge (\ell+s) i \eps-h \eps \} \right\}
\\
   &\ge&  \P \left\{ \bigcup_{i\in I_1} \{ \max_{t\in L_i} Y_1(t) \ge (\ell+s) i \eps-h \eps \} \right\}
   =  \P \left\{ \bigcup_{i\in I_1} \{ X_i^\eps+\delta \xi \ge (\ell+s) i \eps-h \eps \} \right\}.
\end{eqnarray*}
By passing to the complementary events, we obtain
\begin{eqnarray*}
  && \P \left\{ \bigcap_{i\in I_1} \{ X_i^\eps+\delta \xi \le (\ell+s) i \eps-h \eps \} \right\}
  \ge   \P \left\{ \bigcap_{i\in I_1} \{ \tX_i^\eps+\delta \xi_i \le (\ell+s) i \eps-h \eps \} \right\}
\\
  &=& \prod_{i\in I_1}  \P  \left\{ \tX_i^\eps+\delta \xi_i \le (\ell+s) i \eps-h \eps \right\}
  =  \prod_{i\in I_1}  \P  \left\{ X_i^\eps+\delta \xi \le (\ell+s) i \eps-h \eps \right\}.
\end{eqnarray*}
Further, we apply an elementary bound
\begin{eqnarray*}
     1-\P\{\AAA_1\} &=& \P \left\{ \bigcap_{i\in I_1} \{ X_i^\eps \le (\ell+s)i \eps \} \right\}
\\
   &\ge& \P \left\{ \bigcap_{i\in I_1} \{ X_i^\eps+\delta \xi \le (\ell+s)i \eps-h \eps \} \right\}  - \P \{ \delta \xi \le -h \eps \}
\\
    &\ge&  \prod_{i\in I_1}  \P  \left\{ X_i^\eps+\delta \xi \le (\ell+s)i \eps-h \eps \right\} - \P \{ \delta \xi \le -h \eps \}.
\end{eqnarray*}
Under assumption \eqref{h_large}, we have $\P\{\delta \xi\le -h\eps\}\to 0$.
It remains to prove that the product is greater than $\exp(-c\exp(-r))(1+o(1))$.
Taking the logarithm and passing to the complementary events, we see that it is necessary to
prove the bound
\[
   \sum_{i\in I_1} \P \{X_i^\eps+\delta \xi \ge (\ell+s) i \eps-h \eps\}
   \le c\exp(-r) (1+o(1)).
\]
We start with the estimate
\begin{eqnarray} \nonumber
    && \sum_{i\in I_1} \P \{X_i^\eps+\delta \xi \ge (\ell+s) i \eps-h \eps\}
\\    \nonumber
    &\le& \sum_{i\in I_1} \big[ \P \{X_i^\eps \ge (\ell+s) i \eps-2 h \eps\}
    + \P\{\delta \xi>h\eps\}\big]
\\ \label{N1}
 &\le&   \sum_{i: (\ell+s)i+\ell \geq \tau}  \P \{X_i^\eps \ge (\ell+s) i \eps-2 h \eps\}
    + N_1 \P\{\delta \xi> h\eps\},
\end{eqnarray}
where $N_1$ denotes the number of elements in the set $I_1$; it has asymptotics
\[
   N_1\sim \frac{\sigma-\tau}{\ell+s}=\frac{(R-r)B_\eps}{\ell+s}
   \sim \frac{R}{\eps \sqrt{-2\ln \eps}\,\ell}.
\]
For the sum in \eqref{N1} Pickands--Piterbarg lemma provides an equivalent expression
\be \label{sum_pickpit2}
     c \,\ell\, \sum_{i\in I_1} \left( (\ell+s) i \eps-2h\eps \right)^{2/\alpha-1}
     \exp\left(-\left((\ell+s) i \eps -2 h\eps \right)^2/2\right).
\ee
Next, Lemma \ref{l:prop} applied with parameters
$a=\ell+s$, $b=-2h$, $\theta = \tau-\ell-2h$ yields an asymptotics for the latter sum.
Here, as in the derivation of the upper bound we have $a\sim\ell$, $\theta\sim\tau$,
$(\theta\eps)^2=(\tau\eps)^2+o(1)$.
By combining the result of Lemma \ref{l:prop} with \eqref{tauassymp_0}, we obtain
\[
	 \frac{c\,\ell}{a \eps} (\theta\eps)^{2/\alpha-2} \exp\{-(\theta\eps)^2/2\}
     \sim \frac{c}{\eps} (\tau\eps)^{2/\alpha-2}\exp\{-(\tau\eps)^2/2\}
     \sim ce^{-r}.
\]
This means that
\[
    \sum_{i:(\ell+s)i+\ell\ge \tau}  \P \{X_i^\eps \ge (\ell+s) i \eps-2 h \eps\} = c \, e^{-r}(1+o(1)).
\]
It remains to estimate the last term in \eqref{N1}. To this aim, we have to specify the choice of parameters
$\ell$ and $R$.

Since $\P(\delta\xi>h\eps)\to 0$, we may chose $\ell=\ell(\eps)$, although satisfying \eqref{ls2},
but still such that
\[
  \frac{\P(\delta\xi>h\eps)}{\ell \eps\sqrt{-\ln\eps}} \to 0.
\]
Then we may choose $R=R(\eps)$ tending to infinity so slowly that
\[
   N_1 \, \P\{\delta \xi> h\eps\} \sim \frac{ R \, \P(\delta\xi>h\eps)}{\ell \eps \sqrt{-2\ln\eps}} \to 0.
\]
By summing up the estimates for the terms of \eqref{N1}, we arrive at the required lower estimate for $1-\P\{\AAA_1\}$.

\subsection{Proof of Lemma \ref{l:prop}} \label{ss:proof_prop}

The monotone decay of the function $x\mapsto x^{2/\alpha-1}\exp \{- x^2/2\}$ at large $x$
yields the following two bounds
\[
      [(ai+b)\eps]^{2/\alpha-1}\exp\{-[(ai+b)\eps]^2/2\}
      \le \frac{1}{a\eps}  \int_{(a(i-1)+b)\eps}^{(ai+b)\eps}  x^{2/\alpha-1}\exp \{- x^2/2\} dx,
\]

\[
       [(ai+b)\eps]^{2/\alpha-1}\exp\{-[(ai+b)\eps]^2/2\}
       \ge \frac{1}{a\eps}  \int_{(ai+b)\eps}^{(a(i+1)+b)\eps}  x^{2/\alpha-1}\exp \{- x^2/2\} dx.
\]
By summing up over all considered $i$ we infer that the sum under consideration is contained between
two integrals
\[
\frac{1}{a\eps}\int_{(\theta-a)\eps}^\infty x^{2/\alpha-1}\exp \{- x^2/2\} dx \;\; \mbox{ and} \;\;
\frac{1}{a\eps}\int_{(\theta+a)\eps}^\infty x^{2/\alpha-1}\exp \{- x^2/2\} dx.
\]
Furthermore, under our assumptions on $\theta$ and $a$ it is true that
\begin{eqnarray*}
	\frac{1}{a \eps} \int_{(\theta-a) \eps}^\infty x^{2/\alpha-1} \exp \{- x^2/2\} dx &\sim&
	\frac{1}{a \eps}[(\theta-a) \eps]^{2/\alpha-2} \exp\{-[(\theta-a) \eps]^2/2\}
	\\ &\sim& \frac{1}{a \eps}(\theta \eps)^{2/\alpha - 2} \exp\{-(\theta\eps)^2/2\}
\end{eqnarray*}
and
\[
\frac{1}{a\eps}\int_{(\theta+a)\eps}^\infty x^{2/\alpha-1}\exp \{- x^2/2\} dx \sim
\frac{1}{a\eps}(\theta\eps)^{2/\alpha-2} \exp\{-(\theta\eps)^2/2\},
\]
and the required estimate follows.

\end{document}